\documentclass{tac}
\usepackage{xy}
\xyoption{all}

\def\mathcaldef#1{\expandafter\def\csname#1\endcsname{{\cal#1}}}

\def\q{\quad}
\def\qq{\quad\quad}

\def\iso{\,\cong\,}
\def\equ{\,\simeq\,}
\def\inc{\hookrightarrow}
\def\down{\downarrow\!\!}
\def\up{\uparrow\!\!}
\def\la{\langle}
\def\ra{\rangle}
\def\adj{\dashv}
\def\bs{\backslash}
\def\op{^{\rm op}}
\def\ex{\exists}
\def\comp{\pi_0}
\def\ov{\overline}
\def\eps{\varepsilon}

\mathcaldef{C}
\mathcaldef{E}
\mathcaldef{M}
\mathcaldef{S}
\mathcaldef{P}
\mathcaldef{A}
\mathcaldef{B}
\mathcaldef{D}

\mathbfdef{Set}
\mathbfdef{Cat}
\mathbfdef{Gph}
\mathbfdef{Pos}
\mathbfdef{mod}

\mathrmdef{hom}
\mathrmdef{id}
\mathrmdef{ten}

\def\CatX{\Cat\!/\! X}

\def\PrX{\Set^{X\op}}
\def\CX{\C\!/\! X}
\def\CY{\C\!/Y}
\def\MX{\M\!/\! X}
\def\MY{\M\!/Y}
\def\Mu{\M\!/1}
\def\Mou{\M'\!/1}
\def\MoX{\M'\!/\! X}

\def\EM{{\cal(E,M)}}
\def\EMo{{\cal(E',M')}}
\def\Cd{\C^{\rm d}}
\def\Cc{\C^{\rm c}}
\def\SCat{\S\_\Cat^*}
\def\Xmod{X\_\mod}

\newtheorem{prop}{Proposition}
\newtheorem{corol}{Corollary}
\let\pf\proof
\let\epf\endproof
\def\eq{\begin{equation}}
\def\eeq{\end{equation}}

\author{Claudio Pisani}
\address{via Gioberti 86,\\ 10128 Torino, Italy.}
\title{Balanced category theory}
\copyrightyear{2008}
\eaddress{pisclau@yahoo.it}

\begin{document}

\maketitle

\begin{abstract}

Some aspects of basic category theory are developed in a finitely complete category $\C$,
endowed with two factorization systems which determine the same discrete objects and
are linked by a simple reciprocal stability law.
Resting on this axiomatization of final and initial functors and discrete (op)fibrations, 
concepts such as components, slices and coslices, colimits and limits, left and right adjunctible maps, 
dense maps and arrow intervals, can be naturally defined in $\C$, and several classical properties 
concerning them can be effectively proved.

For any object $X$ of $\C$, by restricting $\CX$ to the slices or to the coslices of $X$, two dual 
``underlying categories" are obtained. These can be enriched over internal sets (discrete objects) of $\C$:
internal hom-sets are given by the components of the pullback of the corresponding slice and coslice of $X$.
The construction extends to give functors $\C\to\Cat$, which preserve (or reverse) slices and adjunctible maps
and which can be enriched over internal sets too.
\tableofcontents

\end{abstract}

\section{Introduction}
\label{intro}

Following~\cite{pis3}, we further pursue the goal of an axiomatization of the category of categories
based on the comprehensive factorization system (cfs for short).
Here, the emphasis is on the simultaneous consideration of both the cfs's $\EM$ and $\EMo$ on $\Cat$,
where $\E$ ($\E'$) is the class of final (initial) functors and $\M$ ($\M'$) is the class of discrete
(op)fibrations.
Rather than on the duality functor $\Cat\to \Cat$, we rest on the following ``reciprocal stability" property:
{\em the pullback of an initial (final) functor along a discrete (op)fibration is itself initial (final)}  
(Proposition~\ref{4}).

A ``balanced factorization category" is a finitely complete category $\C$ with two factorization systems 
on it, which are reciprocally stable in the above sense and generate the same subcategory of ``internal sets":
$\S:=\Mu=\Mou\inc\C$.
It turns out that these simple axioms are remarkably powerful, allowing to develop within $\C$ 
some important aspects of category theory in an effective and transparent way.

In particular, the reflection of an object $X\in\C$ in $\S$ gives the internal set of components $\comp X\in\S$,
while the reflections of a ``point" $x:1\to X$ in $\MX$ and $\MoX$ are the slice $X/x$ and the coslice $x\bs X$,
with their maps to $X$.
One can define two functors $\C\to\Cat$, which on $X\in\C$ give the ``underlying" categories
$\ov X(x,y):=\CX(X/x,X/y)$ and $\ov X'(x,y):=\CX(x\bs X,y\bs X)$.
Furthermore, if $X[x,y]:= x\bs X \times_X X/y$, the internal set $X(x,y):=\comp X[x,y]$ 
enriches both $\ov X(x,y)$ and $\ov X'(y,x)$. 
Also composition in $\ov X$ and $\ov X'$ can be enriched to a natural ``symmetrical"
composition map $X(x,y)\times X(y,z)\to X(x,z)$, and it follows in particular that the underlying categories are dual.

In fact, the term ``balanced" is intended to reflect the fact that the objects of $\C$ carry two dual 
underlying categories, coexisting on the same level.
It is the choice of one of the two factorization systems which determines one of them as ``preferred":
in fact, we will make only terminologically such a choice, by calling ``final" (``initial") the maps in $\E$ ($\E'$) 
and ``discrete fibrations" (``opfibrations") those in $\M$ ($\M'$).
For example, we will see that the underlying functor $\C\to\Cat$ determined by $\EM$ preserves discrete fibrations 
and opfibrations, final and initial points, slices and coslices and left and right adjunctible maps, 
while the other one reverses them.

The structure of balanced factorization category also supports the notion of 
``balanced" cylinder and ``balanced" homotopy.
In particular, for any $x,y:1\to X$ we have a cylinder (see ~\cite{law94})
\[ \xymatrix@R=4pc{ X(x,y)\ar@<.5ex>[r]^i\ar@<-.5ex>[r]_e & X[x,y] \ar[r]^c &  X(x,y) } \]
with $i\in\E'$, $e\in\E$ and $c\in\E\cap\E'$,
and pulling back along an element of $\alpha:1\to X(x,y)$ (corresponding to an arrow $x\to y$ as defined above)
we get the ``arrow interval" $[\alpha]$: 
\[ \xymatrix@R=4pc{ 1\ar@<.5ex>[r]^i\ar@<-.5ex>[r]_e & [\alpha] } \]
which in $\Cat$ is the category of factorizations of $\alpha$ with the two factorizations through identities.
 
Another consequence of the axioms is the following exponential law (Proposition~\ref{78}):
if $m\in\MX$, $n\in\MoX$ and the exponential $n^m$ exists in $\CX$, then it is in $\MoX$.
As a consequence, if $m\in\MX$ has a ``complement" (see~\cite{pis} and~\cite{pis2}): 
\[ \neg m : \S \to \CX \qq S\mapsto (S\times X)^m \]
then it is valued in $\MoX$ (and conversely).
In the balanced factorization category of posets (see Section~\ref{ex}), the internal sets are the truth values
and the above enrichment gives the standard identification of posets with $2$-categories,
while the complement operator becomes the usual (classical) one, relating upper-sets and lower-sets.

\section{Some categorical notions}
\label{cat}

We here stress some properties of $\Cat$ and of the comprehensive factorization systems $\EM$ and $\EMo$, 
on which we will model the abstraction of Sections~\ref{emcat} and~\ref{bal}.

\subsection{Slices and intervals}

Recall that given two functors $p:P\to X$ and $q:Q\to X$, the map (or ``comma") construction $(p,q)$  
yields as particular cases the slices $X/x=(\id,x)$ or coslices $x\bs X=(x,\id)$, for any object $x:1\to X$.
More generally, $p/x=(p,x)$ and $x\bs q=(x,q)$ can be obtained as the pullbacks
\eq   \label{1}
\xymatrix@R=4pc@C=4pc{
p/x \ar[r]\ar[d]  & X/x \ar[d] \\
P  \ar[r]^p       & X          }
\qq\qq            
\xymatrix@R=4pc@C=4pc{
x\bs q \ar[r]\ar[d]  & x\bs X \ar[d] \\
Q  \ar[r]^q       & X          }
\eeq
In particular, the pullback
\eq   \label{2}
\xymatrix@R=4pc@C=4pc{
[x,y] \ar[r]\ar[d]  & X/y \ar[d] \\
x\bs X  \ar[r]      & X          }
\eeq
is the category of factorizations between $x$ and $y$, which has as objects consecutive arrows 
$\beta:x\to z$ and $\gamma:z\to y$ in $X$, and as arrows the diagonals $\delta$ as below:
\eq  \label{3}
\xymatrix@R=3.5pc@C=3.5pc{
x \ar[r]^{\beta'}\ar[d]_\beta &  z' \ar[d]^{\gamma'} \\
z \ar[r]^\gamma\ar[ur]^\delta &  y           }
\eeq
Applying the components functor $\comp:\Cat\to\Set$ we get $\comp[x,y]=X(x,y)$.
Indeed $[x,y]$ is the sum $\sum[\alpha], \alpha:x\to y$, where $[\alpha]$ is the ``interval category" of 
factorizations of $\alpha$; $\alpha=\gamma\beta$ is an initial (terminal) object of $[\alpha]$ 
iff $\beta$ ($\gamma$) is an isomorphism.

For any functor $f:X\to Y$ and any $x\in X$, there is an obvious ``slice functor" of $f$ at $x$:
\eq   \label{3b}
e_{f,x}:X/x\to Y/fx
\eeq

\subsection{The comprehensive factorization systems}

Recall that a functor is a discrete fibration if it is orthogonal to the codomain functor $t:1\to 2$,
or, equivalently, if the object mapping of the slice functor $e_{f,x}$ is bijective, for any $x\in X$.
A functor $p:P\to X$ is final if $\comp(x\bs p)=1$ for any $x\in X$.
For example, an object $e:1\to X$ is final iff it is terminal ($\comp(x\bs e)=\comp X(x,e)=\comp 1=1$), 
the identity (or any isomorphism) $X\to X$ is final ($\comp(x\bs\id)=\comp(x\bs X)=1$, since $x\bs X$ has an 
initial object and so it is connected) and $!:X\to 1$ is final iff $X$ is connected ($\comp(0\bs !)=\comp(X)=1$).

Final functors and discrete fibrations form the ``left" comprehensive factorization system $\EM$ 
on $\Cat$, which is in fact the (pre)factorization system generated by $t$. 
Since any $\EM$-factorization gives a reflection in $\MX$ (see Section~\ref{emcat}) we have 
a left adjoint $\down_X\adj i_X$ to the full inclusion of $\MX\equ\PrX$ in $\CatX$.
The ``reflection formula" (see~\cite{pis} and~\cite{pis2}, and the references therein) 
\eq   \label{3a}
(\down p)x \iso \comp(x\bs p) 
\eeq
gives its value at $x$. 
Dually, the domain functor $s:1\to 2$ generates the ``right" cfs $\EMo$: initial functors and discrete opfibrations
(this is the cfs originally considered in~\cite{stw73}).

In particular, the reflection of an object $x:1\to X$ in $\MX$ is the slice projection $\down x:X/x\to X$ 
(corresponding to the representable presheaf $X(-,x)$) and the universal property of the reflection 
reduces to the (discrete fibration version of the) Yoneda Lemma. 
So, a functor to $X$ is isomorphic in $\CatX$ to the slice projection $X/x\to X$ iff it is a discrete fibration
whose domain has a terminal object over $x$.
Dually, if $\up_X\adj j_X:\MoX\to\Cat$ is the reflection in discrete opfibrations, $\up x:x\bs X\to X$ 
in $\MoX$ is the coslice projection.

On the other hand, the reflection of $X\to 1$ in $\Mu$ gives the components of $X$, 
that is, $\down_1 : \Cat/1\iso\Cat \to \Mu \equ \Set$ can be identified with the components functor $\comp$.
The same of course holds for $\up_1 : \Cat \to \Mou \equ \Set$.
\begin{prop}    \label{4}
The pullback of an initial (final) functor along a discrete (op)fibration is itself initial (final).
\end{prop}
\pf
The two statements are clearly dual. 
If in the diagram below $q$ is final, $f$ is a discrete opfibration, both squares are pullbacks
and $x = \down x\circ i$ is a $\EMo$-factorization of $x:1\to X$,
\eq  \label{5}
\xymatrix@R=4pc@C=4pc{
          & x\bs p \ar[r]\ar[d]      & P \ar[r]\ar[d]_p &  Q \ar[d]^q \\
1\ar[r]^i & x\bs X \ar[r]^{\up\,\,x} & X \ar[r]^f       &  Y           }
\eeq
we want to show that $x\bs p$ is connected.
Since $f\circ\up x$ is a discrete opfibration whose domain has an initial object, 
it is a coslice of $Y$: $x\bs X\iso fx\bs Y$ ($fx = (f\circ\!\down x)\circ i$ is a $\EMo$-factorization);
then the vertex of the pullback rectangle is $x\bs p\iso fx\bs q$, which is connected by hypothesis.
\epf

\section{$\EM$-categories and $\EM$ category theory}
\label{emcat}

Throughout the section, $\C$ will be an $\EM$-category, that is a finitely complete category with 
a factorization system on it.

\subsection{Factorization systems and (co)reflections}
\label{fact}

For a brief survey on factorization systems and the associated bifibrations on $\C$
we refer also to~\cite{pis3}.
We just recall that an $\EM$-factorization of a map $f:X\to Y$ in $\C$
\[
\xymatrix@R=3pc@C=3pc{
X  \ar[r]^e\ar[dr]_f  & Z \ar[d]^m \\
                      & Y           }
\]
gives both a reflection of $f\in\CY$ in $\MY$ (with $e$ as reflection map) and a coreflection
of $f\in X\bs\C$ in $X\bs\E$ (with $m$ as coreflection map).
Conversely, any such (co)reflection map gives an $\EM$-factorization.

This can be restated in terms of arrow intervals: say that $f=hg$ is a left (right) point 
of $[f]$ if $g$ ($h$) is in $\E$ ($\M$). 
Then by orthogonality there is a unique map in $[f]$ from any left point to any right one,
and an $\EM$-factorization of $f$ is a ``middle point" of the interval $[f]$: 
it is final among left points and initial among right ones.
 
The above point of view may be reinforced by geometrical intuition (as discussed at length in~\cite{pis3}).
A map $f:X\to Y$ is a figure in $Y$ of shape $X$ (or a cofigure of $X$ in $Y$). 
A left point of $[f]$ factorizes the figure
through an ``infinitesimal modification" of $X$, while a right point factorizes it through
a ``local aspect" of $Y$.   
The middle point of the interval $[f]$ 
\[
\xymatrix@R=3pc@C=3pc{
X  \ar[r]^{e_f} & N(f) \ar[r]^{m_f} & Y }
\]
factorizes the figure through its ``neighborhood" $N(f)$. 
If monic maps are concerned, a left point is an ``infinitesimal enlargement" of $X$
with respect to the embedding $f$ in $Y$, while a right point is an ``open" part of $Y$ which contains $f$.
Then the (infinitesimal) neighborhood $N(f)$ is  both the smallest open part of $Y$
containing $f$ and the biggest infinitesimal enlargement of $X$ with respect to $f$.
In this vein, maps $m:X\to Y$ in $\M$ are characterized as ``local homeomorphisms" 
(if $t:T\to X$ is a figure of $X$, then $m\circ m_t:N(t)\to Y$ gives the neighborhood $N(mt)\to Y$); 
on the other hand, maps $e:X\to Y$ in $\E$ preserve ``global aspects": 
if $h:Y\to Z$ is a cofigure of $Y$ and $e_h$ is the corresponding maximal infinitesimal enlargement of $Y$ in $Z$, 
then $e_h\circ e$ is a maximal ``infinitesimal enlargement" $X\to N(he)$ of $X$ in $Z$.

However, our intuition of ``infinitesimal modification" should be broad enough to include, say, 
the dense inclusion of a part of a topological space in another one, and also
the reflection map of a space in its set of components $X\to\comp X$. 
For example, if $\C=\Cat$, components are one of the ``global aspects" preserved (in the above sense)
by final functors.

\subsection{$\EM$ category theory}
\label{emct}

We now define and analyze several $\EM$-concepts, which assume their usual meaning when $\C=\Cat$ 
with the left comprehensive factorization system  (see also~\cite{pis3}, where the same topics 
are treated in a slightly different way, and other classical theorems are proved).
We say that a map $e\in\E$ is {\bf final}, while a map $m\in\M$ is a {\bf discrete fibration} (df).

We assume that a ``canonical" factorization of any arrow in $\C$ has been fixed, 
and denote by $\down_X\adj i_X$ the corresponding left adjoint to the full inclusion $i_X:\MX\inc\CX$.

\subsection{Sets and components}

If $S\to 1$ is a discrete fibration, we say that $S$ is an {\bf internal set} or also a {\bf $\C$-set} 
(or simply a set) and we denote by $\S := \Mu\inc\C$ the reflective full subcategory of $\C$-sets;
$\S$ is closed with respect to limits in $\C$, and so is itself finitely complete.
The reflection $\comp:=\,\,\down_1:\C\to\S$ is the {\bf components} functor. 
So a final map $e:X\to S$ to a set, that is a factorization of the terminal map $X\to 1$,
gives the set $S \iso \comp X$ of components of $X$, with the following universal property:

\eq   \label{6}
\xymatrix@R=3pc@C=3pc{
X \ar[r]^e\ar[dr]\ar@/^1.5pc/[rr]  & S \ar[d]^m \ar@{..>}[r] & S' \ar[dl]^n \\
                                   & 1                       &               }
\eeq
That is, $\comp X$ is initial among {\em sets} with a map from $X$ and, 
as remarked in Section~\ref{fact}, it is also final among objects with a {\em final} map from $X$.
In particular, a space $X$ is {\bf connected} ($\comp X \iso 1$) iff $X\to 1$  is final, 
iff any map to a set is constant.

\subsection{Slices}

The other important particular case is the reflection $\down x$ of a ``point" (global element) $x:1\to X$,
which is obtained by factorizing the point itself:
\[
\xymatrix@R=3pc@C=3pc{
1  \ar[r]^{e_x}\ar[dr]_x  & X/x \ar[d]^{\down\,\, x} \\
                          & X                     }
\]
We say that $X/x = N(x)$ is the {\bf slice} of $X$ at $x$ (or the ``neighborhood" of $x$, as in 
Section~\ref{fact} and in~\cite{pis3}). 
The corresponding universal property is the following:
\eq   \label{7}
\xymatrix@R=3pc@C=3pc{
1 \ar[r]^{e_x}\ar[dr]_x\ar@/^1.5pc/[rr]^a & X/x \ar[d]^{\down\,\, x} \ar@{..>}[r]^u & A \ar[dl]^m \\
                                          & X                                       &               }
\eeq
which in $\Cat$ becomes the Yoneda Lemma:
given a discrete fibration $m:A\to X$ and an object $a$ of $A$ over $X$, there is a unique factorization $a=ue_x$ 
over $X$ through the ``universal element" $e_x$.
Note that $u$ is a discrete fibration too, and so displays $X/x$ as a slice $A/a$ of $A$.

\subsection{Cones and arrows}

Given a map $p:P\to X$ and a point $x$ of $X$, a {\bf cone} $\gamma:p\to x$ is 
a map over $X$ from $p$ to the slice $X/x$:
\[
\xymatrix@R=3pc@C=3pc{
P  \ar[r]^\gamma\ar[dr]_p  & X/x \ar[d]^{\down\,\, x} \\
                           & X                        }
\]
A cone $\lambda:p\to x$ is {\bf colimiting} if it is universal among cones with domain $p$:
\eq   \label{9}
\xymatrix@R=3pc@C=3pc{
P \ar[r]^\lambda\ar[dr]_p\ar@/^1.5pc/[rr]^\gamma & X/x \ar[d]^{\down\,\, x} \ar@{..>}[r]^u & X/y \ar[dl]^{\down\,\, y} \\
                                                 & X                                       &               }
\eeq
That is, a colimiting cone gives a reflection of $p$ in the full subcategory $\ov X\inc\MX$ 
generated by the slices of $X$ (the {\bf underlying category} of $X$).
If $\gamma:p\to x$ is final, it is an {\bf absolute} colimiting cone.
Note that $\gamma$ is indeed colimiting, since it gives a reflection of $p$ in $\MX$. 

If $P=1$, a cone $\alpha:x\to y$ is an {\bf arrow} from $x$ to $y$:
\[
\xymatrix@R=3pc@C=3pc{
1 \ar[r]^\alpha\ar[dr]_x  & X/y \ar[d]^{\down\,\, y} \\
                          & X                        }
\]
Arrows can be composed via the Kleisli construction: 
if $\alpha:x\to y$, $\beta:y\to z$ and $\widehat\beta$ is universally defined by the left hand diagram below,
then $\beta\circ\alpha$ is given by the right hand diagram:
\eq   \label{10}
\xymatrix@R=3pc@C=3pc{
1 \ar[r]|{e_y}\ar[dr]_y\ar@/^1.5pc/[rr]^\beta & X/y \ar[d]^{\down\,\,y} \ar[r]|{\widehat\beta} & X/z \ar[dl]^{\down\,\,z} \\
                                              & X                                          &               }
\qq\qq   
\xymatrix@R=3pc@C=3pc{
1 \ar[r]|\alpha\ar[dr]_x\ar@/^1.5pc/[rr]^{\beta\,\circ\,\alpha} & X/y \ar[d]^{\down\,\,y} \ar[r]|{\widehat\beta} & X/z \ar[dl]^{\down\,\,z} \\
                                                                & X                                          &               }
\eeq
We so get a category which is clearly isomorphic to the underlying category $\ov X$, and which
we will be identified with it whenever opportune.
Note that the identities $e_x$ are fixed by the choice of a ``canonical" $\EM$-factorization of points.  
An arrow $\alpha:x\to y$ is colimiting iff it is absolutely so, iff it is an isomorphism in $\ov X$.

If $\,\down p:N(p)\to X$ is the reflection of $p:P\to X$ in $\MX$, the reflection map $e_p$ induces  
a bijection between cones $p\to x$ and cones $\down p\to x$, which is easily seen to restrict to 
a bijection between the colimiting ones:
\eq   \label{11}
\xymatrix@R=3pc@C=3pc{
            &                   & X/y\ar@/^.5pc/[ddl]|<<<<<<<<<<{\down\,\,y}       \\
P \ar[r]|{e_p}\ar[dr]_p\ar[rru] & N(p)\ar[ru]\ar[d]^{\down\,\,p} \ar[r]|\lambda  
& X/x \ar@{..>}[u]_u\ar[dl]^{\down\,\,x} \\
            & X                                 &               }
\eeq
\begin{prop}    \label{12} 
If $\,e:T\to P$ is a final map and $p:P\to X$, then a cone $\lambda:p\to x$ is colimiting, 
iff $\,\lambda\circ e:p\circ e\to x$ is colimiting.
\end{prop}
\pf
Since $e_p\circ e: p\circ e \to\,\down p$ is a reflection map of $p\circ e$ in $\MX$,
by the above remark, $\lambda=\lambda'\circ e_p$ is colimiting iff $\lambda'$ is colimiting, 
iff $\lambda'\circ e_p\circ e$ is colimiting.
\eq   \label{13}
\xymatrix@R=3pc@C=3pc{
T\ar[r]|e & P\ar[r]|{e_p}\ar[dr]_p & N(p)\ar[d]^{\down\,\,p} \ar[r]|{\lambda'} & X/x \ar[dl]^{\down\,\,x} \\
          &                        & X                                         &               }
\eeq
\epf

\subsection{Arrow maps}

A map $f:X\to Y$ in $\C$ acts on cones and arrows via an induced map between slices as follows (see~(\ref{3b})).
By factorizing $f\circ\down x$ as in the left hand diagram below, we get a df over $Y$ with a final point
over $fx$. Then there is a (unique) isomorphism $A\to Y/fx$ over $Y$ which respects the selected final points.
By composition, we get the right square below. Its upper edge $e_{f,x}$, which takes the final point of $X/x$ 
to that of $Y/fx$, is the {\bf arrow map} (or {\bf slice map}) of $f$ (at $x$):
\eq   \label{14}
\xymatrix@R=4pc@C=3pc{
1 \ar[r]^{e_x}\ar[dr]_x & X/x \ar[r]^e\ar[d]^{\down\,\,x}  & A \ar[d]^m\ar[r]^\sim & Y/fx \ar[ld]^{\down\,\,fx}   \\
                        & X   \ar[r]^f                     & Y                     &                              }
\qq\qq
\xymatrix@R=4pc@C=4pc{
X/x \ar[r]^{e_{f,x}}\ar[d]^{\down\,\, x} & Y/fx \ar[d]^{\down\,\,fx} \\
X   \ar[r]^f                             & Y                         }
\eeq
Note that $e_{f,x}$ is uniquely characterized, among final maps $X/x\to Y/fx$, by the equation $e_{f,x}\circ e_x = e_{fx}$.
Indeed, given another such $e'_{f,x}$, the dotted arrow induced among the $\EM$-factorizations
is the identity, since it takes the ``universal point" $e_{fx}$ to itself:
\eq    \label{15}
\xymatrix@R=4pc@C=4pc{
&                                         &                                        & Y/fx\ar[ddl]^{\down\,\,fx} \\
1\ar[r]^{e_x} & X/x \ar[r]^{e_{f,x}}\ar[urr]^{e'_{f,x}}\ar[d]^{\down\,\, x} & Y/fx \ar[d]^{\down\,\,fx}\ar@{..>}[ur] & \\
& X   \ar[r]^f                             & Y                                      &      }
\eeq

The arrow map of $f$ acts by composition on cones: it takes $\gamma:p\to x$ to its ``image cone" $f\gamma:fp\to fx$:
\eq   \label{16}
\xymatrix@R=4pc@C=4pc{
P \ar[r]^\gamma\ar[dr]_p\ar@/^1.5pc/[rr]^{f\gamma} & X/x \ar[r]^{e_{f,x}}\ar[d]^{\down\,\, x} & Y/fx \ar[d]^{\down\,\,fx} \\
                                                   & X   \ar[r]^f                             & Y                         }
\eeq
which has the usual meaning in $\Cat$.
\begin{prop}    \label{17}   
An absolute colimiting cone is preserved by any map.
\end{prop}
\pf
If in the diagram above $\gamma$ is final, so is also $f\gamma$.
\epf

\subsection{The underlying functor}

In particular, $f$ induces a mapping between arrows $x\to x'$ of $X$ and arrows $fx\to fx'$ of $Y$, which
preserves identities.
In fact, it is the arrow mapping of a functor $\ov f:\ov X\to\ov Y$: since in the left hand diagram below $e_{f,x}$
is final, there is a uniquely induced dotted arrow $u$ (over $Y$) which makes the upper square commute.
So,
\[ u\circ e_{fx} = u\circ e_{f,x}\circ e_x = e_{f,x'}\circ\widehat\alpha\circ e_x = e_{f,x'}\circ\alpha = f\alpha \]
Thus, $u = \widehat{f\alpha}$ and the right hand square commutes, showing that composition is preserved too. 
\eq     \label{18}
\xymatrix@R=3pc@C=3pc{
X/x \ar@/_1.5pc/[dd] \ar[d]_{\widehat\alpha}\ar[r]^{e_{f,x}} & Y/fx \ar@{..>}[d]^u \ar@/^1.5pc/[dd]   \\
X/x' \ar[d]\ar[r]^{e_{f,x'}}                                 & Y/fx' \ar[d]                         \\
X  \ar[r]^f                                                  & Y                                     }   
\qq\qq
\xymatrix@R=3pc@C=3pc{
X/x \ar[d]_{\widehat\alpha}\ar[r]^{e_{f,x}} & Y/fx \ar[d]^{\widehat{f\alpha}}  \\
X/x'\ar[r]^{e_{f,x'}}                       & Y/fx'                            }
\eeq
Furthermore, $e_{g,fx}\circ e_{f,x} = e_{gf,x}$ (since $e_{g,fx}\circ e_{f,x}\circ e_x = e_{g,fx}\circ e_{fx} = e_{gfx}$),
\eq     \label{18a}
\xymatrix@R=4pc@C=4pc{
X/x \ar[r]^{e_{f,x}}\ar@/^1.5pc/[rr]^{e_{gf,x}}\ar[d]^{\down\,\, x} & Y/fx \ar[r]^{e_{g,fx}}\ar[d]^{\down\,\,fx} & Z/gfx \ar[d]^{\down\,\,gfx}\\
X   \ar[r]^f                                                        & Y    \ar[r]^g                              & Z                          }
\eeq
so that $\ov{g\circ f} = \ov g\circ\ov f$, and we get the {\bf underlying functor} $\ov{(-)}:\C\to\Cat$.
\begin{prop}    \label{18b}   
The underlying functor preserves the terminal object, sets, slices and discrete fibrations.
\end{prop}
\pf
The fact that $\ov 1$ is the terminal category is immediate. 
For any final point $t:1\to X$ in $\C$, the slice projection $X/t\to X$ is an isomorphism,
that is a terminal object in $\CX$ and so also in $\ov X$.
If $f:X\to Y$ is a df in $\C$, its arrow map at any $x$ is an isomorphism:
\eq   \label{19}
\xymatrix@R=4pc@C=4pc{
X/x \ar[r]^\sim \ar[d]^{\down\,\, x}  & Y/fx \ar[d]^{\down\,\,fx} \\
X   \ar[r]^f                          & Y                         }
\eeq
Thus, any arrow $y\to fx$ in $\ov Y$ has a unique lifing along $\ov f$ to an arrow $x'\to x$ in $\ov X$, 
that is $\ov f$ is a discrete fibration in $\Cat$.
The remaining statements now follow at once.
\epf
An example in Section~\ref{ex} shows that connected objects (and so also final maps) are not preserved in general.

\subsection{Adjunctible maps}

If the vertex of the pullback square below has a final point, we say that $f$ is {\bf adjunctible} at $y$:
\eq   \label{20}
\xymatrix@R=3pc@C=3pc{
1 \ar[r]|e\ar[dr]_x\ar@/^1.5pc/[rr]^t & f/y \ar[d]^p \ar[r]|q & Y/y \ar[d]^{\down\,\,y} \\
                                      & X  \ar[r]^f           & Y                       }
\eeq
and that the pair $\la x , t:fx\to y \ra$ is a {\bf universal arrow} from $f$ to $y$.
In that case, $f/y$ is (isomorphic to) a slice $X/x$ of $X$.

Suppose that $f:X\to Y$ is adjunctible at any point of $Y$ and denote by $gy:=pe$ the point of $X$ 
corresponding to (a choice of) the final point $e$, as in the left hand diagram below.
The right hand diagram (whose squares are pullbacks) shows that $g:\C(1,Y)\to\C(1,X)$ can be 
extended to a functor $\,\ov g:\ov Y\to\ov X$. (Note that we are not saying that there is a map $g:Y\to X$ in $\C$
which gives $\ov g$ as its underlying map.) 
\eq   \label{21}
\xymatrix@R=3pc@C=3pc{
1 \ar[r]^{e_{gy}}\ar[dr]_{gy} & X/gy \ar[d]^{\down\,\,gy} \ar[r]^q & Y/y \ar[d]^{\down\,\,y} \\
                              & X   \ar[r]^f                       & Y                       }
\qq\qq   
\xymatrix@R=3pc@C=3pc{
X/gy \ar@/_1.5pc/[dd] \ar@{..>}[d]\ar[r] & Y/y \ar[d] \ar@/^1.5pc/[dd]   \\
X/gy'\ar[d]\ar[r]    & Y/y' \ar[d]                                       \\
X  \ar[r]^f          & Y                                                }
\eeq
\begin{prop}    \label{22}
If $f:X\to Y$ is an adjunctible map, then there is an adjunction $\ov f\adj\ov g:\ov Y\to\ov X$.
\end{prop}
\pf
Suppose that in the left hand diagram below the lower square is a pullback and that $e$ is the arrow map of $f$.
Since $e$ is final, any dotted arrow on the left induces a unique dotted arrow on the right
which makes the upper square (over $Y$) commute.
Since the lower square is a pullback, the converse also holds true.
The bijection is natural, as shown by the diagrams on the right.
\eq   \label{23}
\xymatrix@R=3pc@C=3pc{
X/x \ar@/_1.5pc/[dd] \ar@{..>}[d]\ar[r]^e & Y/fx \ar@{..>}[d] \ar@/^1.5pc/[dd]   \\
X/gy\ar[d]\ar[r]                          & Y/y \ar[d]                           \\
X  \ar[r]^f                               & Y                                     }
\qq\qq 
\xymatrix@R=3pc@C=3pc{
X/x \ar@/_1.5pc/[dd] \ar@{..>}[d]\ar[r]^e & Y/fx \ar@{..>}[d] \ar@/^1.5pc/[dd]   \\
X/gy\ar[d]\ar[r]                          & Y/y \ar[d]                           \\
X/gy'  \ar[r]^f                           & Y/y'                                  }
\qq\qq 
\xymatrix@R=3pc@C=3pc{
X/x' \ar@/_1.5pc/[dd] \ar[d]\ar[r]^{e'} & Y/fx' \ar[d] \ar@/^1.5pc/[dd]   \\
X/x  \ar@{..>}[d]\ar[r]^e               & Y/fx \ar@{..>}[d]               \\
X/gy  \ar[r]^f                          & Y/y                              }
\eeq
\epf
\begin{corollary}    \label{24}
The underlying functor preserves adjunctible maps.
\end{corollary}
\epf
\begin{remarks}
\begin{enumerate}
\item
Diagram~(\ref{20}) shows that Corollary~\ref{24} would be obvious if, along with slices, 
the underlying functor would preserve also pullbacks.
Up to now, the author does not know if this is true in general.
\item
The adjunction of Proposition~\ref{22} is nothing but the adjunction $\ex_f\adj f^*$ 
(where $\ex_f:\MX\to\MY$ gives, when $\C=\Cat$, the left Kan extension along $f$), 
which for an adjunctible $f$ restricts to $\ov f \adj \ov g$ (see~\cite{pis3}).
\item
The diagrams below display the unit and the counit of adjunction $\ov f\adj\ov g$:   
\eq    \label{25}
\xymatrix@R=3pc@C=3pc{
X/gfx \ar[drr]^q\ar[ddr]_{\down\,\,gfx} \\
           &  X/x \ar@{..>}[ul]|-u\ar[r]|{e_{f,x}}\ar[d]^{\down\,\, x} & Y/fx \ar[d]^{\down\,\,fx} \\
           &  X   \ar[r]^f                                             & Y                         }
\qq
\xymatrix@R=3pc@C=3pc{
                                                  &  &  Y/fgy \ar@{..>}[dl]|-v\ar[ddl]^{\down\,\,fgy} \\
X/gy \ar[r]|q\ar[d]^{\down\,\,gy}\ar[urr]^{e_{f,gy}} & Y/y \ar[d]^{\down\,\,y}                        \\
X   \ar[r]^f                                         & Y                                              }
\eeq
In $\Cat$, $u$ is the functor $\alpha:x'\to x \q\mapsto\q \eta_x\circ\alpha:x'\to gfx$, 
and the triangle $e_{x,y} = q\circ u$ reduces to the classical expression of the arrow mapping
of a left adjoint functor via transposition and the unit $\eta$.
Similarly, $v$ is the functor $\beta \mapsto \eps_y\circ\beta$, 
and the triangle $q = v\circ e_{f,gy}$ becomes the classical expression of the transposition bijection
of an adjunction via the arrow mapping of the left adjoint functor and the counit $\eps$.
\end{enumerate}
\end{remarks}
\begin{prop}    \label{26}
An adjunctible map $f:X\to Y$ preserves colimits: 
if $\lambda:p\to x$ is colimiting, then so is $f\lambda:fp\to fx$.
\end{prop}
\pf
Suppose that $\lambda:p\to x$ is colimiting.
We want to show that $f\lambda$ is colimiting as well, that is any cone $\gamma:fp\to y$ 
factorizes uniquely through $e_{f,x}\circ\lambda:fp\to fx$.
\eq   \label{27}
\xymatrix@R=1.5pc@C=1.5pc{
& P \ar[dr]^{f\lambda}\ar@{..>}[dd]|<<<<<u\ar@/^1.5pc/[dddd]|<<<<<<<<<<<<<p\ar[dl]_\lambda\ar@/^2pc/[ddrr]^\gamma  \\
X/x \ar@{..>}[dr]^{u'}\ar[rr]^<<<<<<<{e_{f,x}}\ar[dddr]  & & Y/fx \ar[dddr]\ar@{..>}[dr]^{u''}  \\
& X/gy \ar[rr]^l\ar[dd]        & & Y/y \ar[dd]                \\ \\
& X  \ar[rr]^f                 & & y                            }
\eeq
Since the rectangle is a pullback, we have a universally induced cone $u:p\to gy$, and since $\lambda$ 
is colimiting also $u':X/x\to X/gy$, which by adjunction ($e_{f,x}$ is final) corresponds to a 
unique $u'':Y/fx\to Y/y$, which is the desired unique factorization of $\gamma$ through $f\lambda$.
Note that we have followed the same steps of the most obvious classical proof,
but at this abstract level we need not to verify or use any naturality condition.  
\epf

\subsection{Dense and fully faithful maps}

If in the pullback below $e$ is colimiting we say that $f$ is {\bf adequate} at $y$:
\eq   \label{28}
\xymatrix@R=3pc@C=3pc{
 f/y \ar[d]\ar[r]^e & Y/y \ar[d]^{\down\,\,y} \\
 X  \ar[r]^f        & Y                       }
\eeq
while if $e$ is final we say that $f$ {\bf dense}.
Note that in that case $e$ is an absolute colimiting cone
and the diagram above displays at once dense maps as ``absolutely adequate" maps, 
as well as ``locally final" maps (see~\cite{dense} for the case $\C=\Cat$, where
a different terminology is used; see also Section~\ref{dense}). 

The density condition can be restated as the fact that the counit of the 
adjunction $f^*\adj\ex_f:\MX\to\MY$, that is the uniquely induced arrow $v$ 
on the left, is an isomorphism on the slices of $Y$:
\eq     \label{30}
\xymatrix@R=3pc@C=3pc{
                                &  & \ex_f(f/y) \ar@{..>}[dl]|-v\ar[ddl]^{m} \\
f/y \ar[r]|q\ar[d]^p\ar[urr]^e  & Y/y \ar[d]^{\down\,\,y}                   \\
X   \ar[r]^f                    & Y                                         }
\qq
\xymatrix@R=3pc@C=3pc{
f/fx \ar[drr]^q\ar[ddr]_p                                                                 \\
           &  X/x \ar@{..>}[ul]|-u\ar[r]|{e_{f,x}}\ar[d]^{\down\,\, x} & Y/fx \ar[d]^{\down\,\,fx} \\
           &  X   \ar[r]^f                                             & Y                         }
\eeq
On the other hand, if the unit $u$ of the same adjunction, that is the uniquely induced arrow $u$ 
on the right, is an isomorphism on the slice $X/x$, we say that $f$ is {\bf fully faithful} at $x$.
Note that if $f$ is fully faithful at $x$, then $f$ is adjunctible at $fx$.

\section{Balanced factorization categories}
\label{bal}

As shown in the previous section and in~\cite{pis3}, category theory can in part be developed 
in any $\EM$-category.
However, in order to treat ``left-sided" and ``right-sided" categorical concepts simultaneously
and to analyze their interplay, we need to consider {\em two} factorization systems, satisfying appropriate axioms:
\begin{definition}
A {\bf balanced factorization category} {\rm (bfc)} is a finitely complete category $\C$ with two 
factorization systems $\EM$ and $\EMo$ on it, such that $\Mu=\Mou\inc\C$ and satisfying the 
{\bf reciprocal stability law} {\rm (rsl)}: \\
the pullback of a map $e\in\E$ ($e'\in\E'$) along a map $m'\in\M'$ ($m\in\M$) is itself in $\E$ ($\E'$).
\end{definition}
While this definition aims at capturing some relevant features of $\Cat$, other intersesting
instances of bfc are presented in Setion~\ref{ex}.    

Note that any slice $\CX$ of a bfc is itself a ``weak" bfc: the requirement $\Mu=\Mou$ may not hold.

\subsection{Notations and terminology}

As in Section~\ref{emcat}, we assume that ``canonical" $\EM$ and $\EMo$-factorizations 
have been fixed for the arrows of $\C$, 
and denote by $\down_X:\CX\to\MX$ and $\up_X:\CX\to\MoX$ the corresponding reflections. 
For the concepts defined in Section~\ref{emcat}, but referred to $\EMo$, we use the standard ``dual"
terminology when it is available.
Otherwise, we distinguish the $\EM$-concepts from the $\EMo$-ones, by qualifying them with the
attributes ``left" and ``right" respectively.

So, we say that a map $i\in\E'$ is {\bf initial} and a map $n\in\M'$ is a 
{\bf discrete opfibration} (dof).
 The {\bf coslice} of $X\in\C$ at $x:1\to X$ is:
\[
\xymatrix@R=3pc@C=3pc{
1  \ar[r]^{i_x}\ar[dr]_x  & x\bs X \ar[d]^{\up\,\, x} \\
                          & X                         }
\]
We denote by $\B:=\M\cap\M'$ the class of ``bifibrations".  
Of course, we also have a ``right" underlying functor $\ov{(-)}':\C\to\Cat$,
with $\ov X'(x,y):=\CX(x\bs X,y\bs X)$.

\subsection{First properties}

We begin by noting that, since $\S:=\Mu=\Mou\inc\C$ (that is, ``left" and ``right" $\C$-sets coincide,
and are bifibrations over $1$), then also any (op)fibration over a set is in fact a bifibration:
$\M/S=\M'/S=\B/S=\S/S$. 
Furthermore, a map $X\to S$ to a set is initial iff it is final: $\E/S=\E'/S$.
In that case, $S\iso\comp X$, where $\comp:=\,\up_1 = \,\down_1 :\C\to\S$ is the component functor.
In particular, $\E/1=\E'/1\inc\C$ is the full subcategory of {\bf connected} objects.
Note also that a map $S\to S'$ between sets is final (or initial) iff it is an isomorphism; 
so, if $e:S\to X$ is a final (or initial) map from a set, then $S\iso\comp X$.

\subsection{Cylinders and intervals}

If $c:X\to S$ is a reflection in $\S$, that is if $c$ is initial (or equivalently final),
for any element $s:1\to S$ (which is in $\B$), the pullback
\eq   \label{31}
\xymatrix@R=4pc@C=4pc{
[s]  \ar[r]\ar[d]    & X \ar[d]     \\
1\ar[r]^s            & S               }
\eeq
displays $[s]$ as the ``component" $[s]\inc X$ of $X$ at $s$; $[s]$ is connected by the rsl
and is included in $X$ as a discrete bifibration.
(Note that we are not saying that $X$ is the sum of its components, although this is the case under
appropriate hypothesis.)
If furthermore $i:S'\to X$ is initial (or final), then $S'\iso S\iso\comp X$.
Thus, $i$ is monomorphic and meets any component $[s]$ in a point, which is an initial point $i_s$ of $[s]$:
\eq   \label{32}
\xymatrix@R=4pc@C=4pc{
1\ar[r]^{i_s}\ar[d] & [s] \ar[r]\ar[d]   & 1 \ar[d]^s  \\
S'\ar[r]^i          & X   \ar[r]         & S           }
\eeq
By a {\bf balanced cylinder}, we mean a cylinder 
\[ \xymatrix@R=4pc{ B\ar@<.5ex>[r]^i\ar@<-.5ex>[r]_e & X \ar[r]^c &  B } \]
(that is, a map with two sections; see~\cite{law94}) 
where the base inclusions $i$ and $e$ are initial and final respectively.
As a consequence, $c$ is both initial and final, and if the base is a $\C$-set $S$, 
then $S\iso\comp X$ and $c$ is a reflection in $\S$.
Among these ``discrete cylinders", there are ``intervals", that is balanced cylinders with a terminal base.
Thus, intervals in $\C$ are (connected) objects $i,e:1\to I$ with an initial and a final selected point.

In particular, any component of a cylinder with discrete base is an interval:
\eq   \label{33}
\xymatrix@R=4pc@C=4pc{
1\ar@<.5ex>[r]^{i_s}\ar@<-.5ex>[r]_{e_s}\ar[d] & [s] \ar[r]\ar[d]   & 1 \ar[d]^s \\
S\ar@<.5ex>[r]^i\ar@<-.5ex>[r]_e              & X         \ar[r]   & S          }
\eeq
(Conversely, if $\E$ and $\E'$ are closed with respect to products in $\C^\to$, by multiplying
any object $B\in\C$ with an interval, we get a cylinder with base $B$.)

\subsection{Homotopic maps}

Two parallel maps $f,f':X\to Y$ are {\bf homotopic} if there is a balanced cylinder 
\[ \xymatrix@R=4pc{ X\ar@<.5ex>[r]^i\ar@<-.5ex>[r]_e & C \ar[r]^c &  X } \]
and a map (\v homotopy") $h:C\to Y$ such that $f=hi$ and $f'=he$ (see~\cite{law94}). 
\begin{prop}    \label{34}
The two base inclusions $i,e:X\to C$ of a balanced cylinder are coequalized by any map $C\to S$ toward a set. 
\end{prop}
\pf
Since the retraction $c$ of the cylinder on its base is final,
the proposition follows immediately from the fact that, by orthogonality, any map $C\to S$ to a set
factors uniquely through any final map $c:C\to X$:
\eq    \label{35}
\xymatrix@R=4pc@C=4pc{
C \ar[r]\ar[d]_c       & S \ar[d] \\
X \ar[r]\ar@{..>}[ur]  & 1         }
\eeq
\epf 
\begin{corol}     \label{36}
Two homotopic maps $f,f':X\to Y$ are coequalized by any map $Y\to S$ toward a set. 
\end{corol}
\epf

\subsection{Internal hom-sets}

As in Section~\ref{cat}, we denote by $X[x,y]$ (or simply $[x,y]$) the vertex of a (fixed) pullback
\eq   \label{37}
\xymatrix@R=4pc@C=4pc{
[x,y] \ar[r]^{q_{x,y}}\ar[d]_{p_{x,y}}  & X/y \ar[d]^{\down\,\,y} \\
x\bs X  \ar[r]^{\up\,\,x}               & X          }
\eeq
We also denote by $c_{x,y}:X[x,y]\to X(x,y)$ the ``canonical" reflection map in $\S$.
Thus $X(x,y)$ (or simply $(x,y)$) is the $\C$-set $\comp X[x,y]$.

By further pulling back along the final (initial) point of the (co)slice, we get:
\eq   \label{38}
\xymatrix@R=4pc@C=4pc{
                & S'\ar[r]\ar[d]^e   & 1 \ar[d]^{e_y} \\
S\ar[r]^i\ar[d] & [x,y] \ar[r]\ar[d] & X/y \ar[d]     \\
1\ar[r]^{i_x}   & x\bs X  \ar[r]     & X              }
\eeq
where $S$ and $S'$ are sets (since the pullback of a discrete (op)fibration is still a discrete (op)fibration)
and $i$ ($e$) is initial (final) by the rsl.
Then $b_i=c_{x,y}\circ i:S\to (x,y)$ and $b_e=c_{x,y}\circ e:S'\to (x,y)$ are both isomorphisms 
and in the above diagram we can replace $S$ and $S'$ by $(x,y)$, $i$ by $i\circ b_i^{-1}$ and $e$ 
by $e\circ b_e^{-1}$. 
So we have the following diagram of pullback squares:
\eq   \label{39}
\xymatrix@R=4pc@C=4pc{
                            & (x,y)\ar[r]\ar[d]^{e_{x,y}}            & 1 \ar[d]^{e_y}         \\
(x,y)\ar[r]^{i_{x,y}}\ar[d] & [x,y] \ar[r]^{q_{x,y}}\ar[d]^{p_{x,y}} & X/y \ar[d]^{\down\,\,y} \\
1\ar[r]^{i_x}               & x\bs X  \ar[r]^{\up\,\,x}              & X                       }
\eeq
and the balanced cylinder:
\eq   \label{40}
\xymatrix@C=3pc{ (x,y)\ar@<.5ex>[r]^{i_{x,y}}\ar@<-.5ex>[r]_{e_{x,y}} & [x,y] \ar[r]^{c_{x,y}} &  (x,y) }                    
\eeq
Considering in particular the pullbacks
\eq   \label{41}
\xymatrix@R=4pc@C=4pc{
(x,y)\ar[r]^{q_{x,y}\circ\,i_{x,y}}\ar[d] & X/y \ar[d]     \\
1\ar[r]^x                                 & X               }
\qq
\xymatrix@R=4pc@C=4pc{
(x,y)\ar[r]^{p_{x,y}\circ\,e_{x,y}}\ar[d] & x\bs X \ar[d]  \\
1\ar[r]^y            & X                                   }
\eeq
we see that the elements of $\alpha:1\to X(x,y)$ correspond bijectively to ``left" 
arrows $x\to y$ and to ``right" arrows $y\to x$: 
\eq   \label{41b}
\C(1,X(x,y)) \iso \ov X(x,y) \iso {\ov X}'(y,x)  
\eeq
Then $X(x,y)$ deserves to be called the ``internal hom-set of arrows from $x$ to $y$". 
(Note that by using $[x,y]$ and $(x,y)$ we have notationally made an arbitrary choice among one of two 
possible orders, causing an apperent asymmetry in the subsequent theory.)
From now on, also an element $\alpha:1\to X(x,y)$ will be called an ``arrow from $x$ to $y$". 
Composing it with $q_{x,y}\circ i_{x,y}$ ($p_{x,y}\circ e_{x,y}$) we get the corresponding left (right) arrow
(which will be sometimes denoted with the same name).

As in~(\ref{33}), to any arrow $\alpha:1\to X(x,y)$ there also corresponds an ``arrow interval" $[\alpha]$,
the component at $\alpha$ of the balanced cylinder~(\ref{40}):
\eq   \label{42}
\xymatrix@R=4pc@C=4pc{
1\ar@<.5ex>[r]^i\ar@<-.5ex>[r]_e\ar[d]^\alpha        & [\alpha] \ar[r]\ar[d]   & 1 \ar[d]     \\
(x,y)\ar@<.5ex>[r]^{i_{x,y}}\ar@<-.5ex>[r]_{e_{x,y}} & [x,y]          \ar[r]   & (x,y)         }
\eeq
\begin{prop}    \label{43}
Two points $x,y:1\to X$ are homotopic iff there is an arrow from $x$ to $y$.
\end{prop}
\pf
If $\alpha:x\to y$ is an arrow of $X$, then $[\alpha]\to X$ gives a homotopy from $x$ to $y$.
Conversely, given a homotopy
\( \xymatrix@R=4pc{ 1\ar@<.5ex>[r]^i\ar@<-.5ex>[r]_e & C \ar[r]^h & X } \)
there is a (unique) arrow $i\to e$ of $C$, and so its image under $h$ is the desired arrow of $X$.
\epf

\subsection{Enriched composition}
\label{enrich}

The ``enriched composition" map 
\[ \mu_{x,y,z}:X(x,y)\times X(y,z)\to X(x,z) \] 
is defined as $\mu_{x,y,z}:=c_{x,z}\circ\mu'_{x,y,z}$, where $\mu'$ is universally induced by the pullback:
\eq   \label{44}
\xymatrix@R=2pc@C=2pc{
& (x,y)\times (y,z) \ar[ddl]_l\ar[ddr]^r\ar@{..>}[d]^{\mu'} &   \\
& [x,z] \ar[dl]\ar[dr]\ar[d]^{c_{x,z}} &  \\
x\bs X \ar[dr] & (x,z) & X/z \ar[dl]      \\
& X &                              }
\eeq
in which the maps $l$ and $r$ are the product projections followed by $p_{x,y}\circ e_{x,y}$ 
and $q_{y,z}\circ i_{y,z}$, respectively.
\begin{prop}    \label{44a}
The map $\mu_{x,y,z}$ induces on the elements the composition mapping $\ov X(x,y)\times \ov X(y,z)\to \ov X(x,z)$ 
in $\ov X$.
\end{prop}
\pf
For any arrow $\beta:1\to X(y,z)$ we have a map 
\( [x,\beta]:[x,y]\to[x,z] \): 
\eq   \label{45}
\xymatrix@R=2pc@C=2pc{
[x,y] \ar@{..>}[dr]^{[x,\beta]}\ar[rr]^{q_{x,y}}\ar[dddr]_{p_{x,y}} & & X/y \ar[dddr]\ar[dr]^{\widehat\beta} \\
& [x,z] \ar[rr]^{q_{x,z}}\ar[dd]^{p_{x,z}} & & X/z \ar[dd]     \\ \\
& x\bs X  \ar[rr]    & & X          }
\eeq
In the diagram below, the left dotted arrow (induced by the universality of the lower pullback square) 
and the right one ($(x,\beta):=\comp[x,\beta]:(x,y)\to (x,z)$, induced by the universality of the 
reflection $c_{x,y}:[x,y]\to (x,y)$) are the same, since the horizontal edges of the rectangle are identities:
\eq   \label{46}
\xymatrix@R=4pc@C=4pc{
(x,y) \ar@/_1.5pc/[dd]\ar[r]^{i_{x,y}}\ar@{..>}[d]^{(x,\beta)} & [x,y]\ar[r]^{c_{x,y}}\ar[d]^{[x,\beta]} & (x,y) \ar@{..>}[d]^{(x,\beta)} \\
(x,z) \ar[r]^{i_{x,z}}\ar[d]                                   & [x,z]\ar[r]^{c_{x,z}}\ar[d]^{p_{x,z}}   & (x,z)                \\
1     \ar[r]^{i_x}                                             & x\bs X                                              }
\eeq
So we also have the commutative diagram
\eq   \label{47}
\xymatrix@R=4pc@C=4pc{
1 \ar[r]^\alpha & (x,y) \ar[r]^{i_{x,y}}\ar[d]^{(x,\beta)} & [x,y]\ar[r]^{q_{x,y}}\ar[d]^{[x,\beta]} &  X/y \ar[d]^{\widehat\beta}  \\
                & (x,z) \ar[r]^{i_{x,z}}                   & [x,z]\ar[r]^{q_{x,z}}                   &  X/z                     }   
\eeq
showing that $(x,\beta)$ acts on an arrow $\alpha:1\to X(x,y)$ in the same way as $\widehat\beta$
acts on the corresponding left arrow $q_{x,y}\circ i_{x,y}\circ\alpha:1\to X/y$, that is as $-\circ\beta:\ov X(x,y)\to\ov X(x,z)$.

Now observe that since $(x,\beta)=c_{x,z}\circ[x,\beta]\circ i_{x,y}$, and since $i_{x,y}$ 
and $e_{x,y}$ are trivially homotopic, also $(x,\beta)=c_{x,z}\circ[x,\beta]\circ e_{x,y}$.
But composing $[x,\beta]\circ e_{x,y}:(x,y)\to[x,z]$ with the pullback projections $p_{x,z}$ and $q_{x,z}$
we get $p_{x,y}\circ e_{x,y}$ and the constant map through $\beta:1\to X/y$, so that
the map $[x,\beta]\circ e_{x,y}$ coincides with the composite
\[ \xymatrix@C=4pc{ (x,y) \ar[r]^-\sim & (x,y)\times 1 \ar[r]^-{\id\times\beta} & (x,y)\times(y,z)\ar[r]^-{\mu'} & [x,z]} \]
Summing up, we have seen that:
\[ (x,\beta) = \comp[x,\beta] = c_{x,z}\circ[x,\beta]\circ i_{x,y} = c_{x,z}\circ[x,\beta]\circ e_{x,y} 
= \mu\circ (\id\times\beta) \circ \la \id,!\,\ra \]
So, $\mu(\alpha,\beta) = (x,\beta)\circ\alpha$ corresponds to $\ov\mu(\alpha,\beta)$,
where $\ov\mu$ is the composition in $\ov X$, proving the thesis.
\epf
``Dually", for any $\alpha:x\to y$ there are maps $[\alpha,z]:[y,z]\to[x,z]$ 
and $(\alpha,z):(y,z)\to(x,z)$ such that
\[ (\alpha,z) = \comp[\alpha,z] = c_{x,z}\circ[\alpha,z]\circ i_{y,z} = c_{x,z}\circ[\alpha,z]\circ e_{y,z}
= \mu\circ (\alpha\times\id) \circ \la !,\id \,\ra \]
So, $\mu(\alpha,\beta) = (\alpha,z)\circ\beta$ corresponds to ${\ov\mu}'(\beta,\alpha)$,
where ${\ov\mu}'$ is the composition in ${\ov X}'$.
\begin{corollary}    \label{48}
The two underlying categories are related by duality: ${\ov X}'\iso {\ov X}\op$.
\end{corollary}
\epf

\subsection{Internally enriched categories?}
\label{weak}

Up to now, the author has not been able to prove the associativity of the enriched composition,
nor to find a counter-example.
On the other hand, that the identity arrow $u_x:=c_{x,x}\circ u'_x : 1\to (x,x)$
\eq   \label{49}
\xymatrix@R=2pc@C=3pc{
& 1 \ar[ddl]_{i_x}\ar[ddr]^{e_x}\ar@{..>}[d]^{u'_x} &  \\
& [x,x] \ar[dl]\ar[dr]\ar[d]^{c_{x,x}}              &  \\
x\bs X \ar[dr] & (x,x) & X/x \ar[dl]                   \\
& X &                                                   }
\eeq
satisfies the unity laws for the enriched composition, follows from the previous section
or can be easily checked directly.

\subsection{The action of a map on internal hom-sets}

Given a map $f:X\to Y$ in $\C$, we have an induced map $f_{[x,y]}:X[x,y]\to Y[fx,fy]$:
\eq   \label{50}
\xymatrix@R=2pc@C=2pc{
[x,y] \ar@{..>}[dr]^{f_{[x,y]}}\ar[rr]^{q_{x,y}}\ar[dd]_{p_{x,y}} & & X/y \ar[dd]\ar[dr]^{e_{f,y}} \\
& [fx,fy] \ar[rr]^<<<<<<<<<<{q_{fx,fy}}\ar[dd]|<<<<<<<{p_{fx,fy}} & & Y/fy \ar[dd]                  \\
x\bs X \ar[dr]_{i_{f,x}}  \ar[rr]                                 & & X \ar[dr]^f                    \\
& fx\bs Y  \ar[rr]                                                & & Y                              }
\eeq
In the diagram below the left dotted arrow (induced by the universality of the lower pullback square) 
and the right one ($f_{x,y}:=\comp f_{[x,y]}:X(x,y)\to Y(fx,fy)$, induced by the universality of the 
reflection $c_{x,y}:[x,y]\to (x,y)$) are the same, since the horizontal edges of the rectangle are identities.
\eq   \label{51}
\xymatrix@R=4pc@C=4pc{
(x,y) \ar@/_2pc/[dd]\ar[r]^{i_{x,y}}\ar@{..>}[d]^{f_{x,y}} & [x,y]\ar[r]^{c_{x,y}}\ar[d]^{f_{[x,y]}} & (x,y) \ar@{..>}[d]^{f_{x,y}} \\
(fx,fy) \ar[r]^{i_{fx,fy}}\ar[d]                           & [fx,fy]\ar[r]^{c_{fx,fy}}\ar[d]^{p_{fx,fy}}   & (fx,fy)   \\
1     \ar[r]^{i_{fx}}                                      & fx\bs Y                                                    }
\eeq
So we also have the commutative diagram
\eq   \label{52}
\xymatrix@R=4pc@C=4pc{
1 \ar[r]^\alpha & (x,y) \ar[r]^{i_{x,y}}\ar[d]^{f_{x,y}} & [x,y]\ar[r]^{q_{x,y}}\ar[d]^{f_{[x,y]}} &  X/y \ar[d]^{e_{f,y}}  \\
                & (fx,fy) \ar[r]^{i_{fx,fy}}             & [fx,fy]\ar[r]^{q_{fx,fy}}               &  Y/fy                    }   
\eeq
showing that $f_{x,y}$ acts on an arrow $\alpha:1\to X(x,y)$ in the same way as the arrow map of $f$ 
acts on the corresponding left arrow $q_{x,y}\circ i_{x,y}\circ\alpha:1\to X/y$.
So, $f_{x,y}:X(x,y)\to Y(fx,fy)$ enriches $\ov f_{x,y}:\ov X(x,y)\to \ov Y(\ov fx,\ov fy)$.
We also have the ``dual" commuting diagrams, with $e_{x,y}$ in place of $i_{x,y}$ and so on,
showing that $f_{x,y}$ also enriches $\ov f'_{y,x}:\ov X'(y,x)\to \ov Y'(\ov fy,\ov fx)$.
\begin{corollary}    \label{53}
The two underlying functors are related by duality: $\ov{(-)}'\iso\ov{(-)}\op$.
\end{corollary}
\begin{prop}     \label{55}
The maps $f_{x,y}$ are ``functorial" with respect to enriched composition.
\end{prop}
\pf
We want to show that the following diagram commutes:
\eq   \label{56}
\xymatrix@R=4pc@C=6pc{
(x,y)\times(y,z) \ar[r]^-{f_{x,y}\times f_{y,z}}\ar[d]^\mu & (fx,fy)\times(fy,fz) \ar[d]^\mu  \\
(x,z) \ar[r]^{f_{x,z}}                                     &       (fx,fz)                     }   
\eeq
Since the lower rectangle below commutes
\eq   \label{57}
\xymatrix@R=4pc@C=6pc{
(x,y)\times(y,z) \ar[r]^-{f_{x,y}\times f_{y,z}}\ar[d]^{\mu'} &  (fx,fy)\times(fy,fz) \ar[d]^{\mu'}  \\
[x,z] \ar[r]^{f_{[x,z]}}\ar[d]^c                              &  [fx,fz]\ar[d]^c                     \\
(x,z) \ar[r]^{f_{x,z}}                                        &  (fx,fz)                             }   
\eeq
we have to show that the upper one also commutes,
that is that composing with the pullback projections (say, the second one $q_{fx,fz}:[fx,fz]\to Y/fz$)
we get the same maps:
\eq   \label{58}
\xymatrix@R=4pc@C=4pc{
(x,y)\times(y,z) \ar[r]^-\pi\ar[d]^{\mu'} &  (y,z) \ar[d]|{q_{y,z}\circ\,\,i_{y,z}}  \\
[x,z] \ar[r]^{q_{x,z}}\ar[d]^{f_{[x,z]}}  &  X/z \ar[d]^{e_{f,z}} \\
[fx,fz] \ar[r]^{q_{fx,fz}}                &  Y/fz                  }   
\qq\qq
\xymatrix@R=4pc@C=4pc{
(x,y)\times(y,z) \ar[r]^-\pi\ar[d]|{f_{x,y}\times f_{y,z}}  &  (y,z) \ar[d]^{f_{y,z}}  \\
(fx,fy)\times(fy,fz) \ar[r]^-\pi\ar[d]^{\mu'}               & (fy,fz) \ar[d]|{q_{fy,fz}\circ\,\,i_{fy,fz}} \\
[fx,fz] \ar[r]^{q_{fx,fz}}                                  &  Y/fz                     }   
\eeq
Indeed, this is the case because all the squares commute and 
the composites on the right are the two paths of the rectangle~(\ref{52}).
\epf

\subsection{Enriching discrete fibrations}

To any discrete fibration $\,m:A\to X$ in $\C$ there corresponds a presheaf $\ov m$ on $\ov X$,
as can be seen in several ways:
\begin{enumerate}
\item
As we have seen in Proposition~\ref{18b}), by applying $\ov{(-)}:\C\to\Cat$ to $\,m$, 
we obtain a discrete fibration $\,\ov m:\ov A\to\ov X$ in $\Cat$, corresponding to a presheaf $\ov m:\ov X\op\to\Set$.
\item
Via the inclusion $i:\ov X\inc\MX$, any object $A\in\MX$ is ``interpreted" as a presheaf $\MX(i-,A)$ on $\ov X$.
\item
To any object $x:1\to X$ of $\ov X$ there corresponds the set $\CX(x,m)$ of points of $A$ over $x$,
and to any arrow $\alpha:x\to y$ in $\ov X$ there corresponds a mapping $\ov m\alpha:\CX(y,m)\to\CX(x,m)$, 
whose value at $a\in A$ over $y$ is given by the right hand diagram below
(where $\widehat a$ is defined by the left hand one):  
\eq     \label{59}
\xymatrix@R=3pc@C=3pc{
1 \ar[r]^{e_y}\ar[dr]_y\ar@/^1.5pc/[rr]^a & X/y \ar[d]^{\down\,\,y} \ar[r]^{\widehat a} & A \ar[dl]^m \\
                                          & X                                           &               }
\qq\qq   
\xymatrix@R=3pc@C=3pc{
1 \ar[r]^\alpha\ar[dr]_x\ar@/^1.5pc/[rr]^{(\ov m\alpha)a} & X/y \ar[d]^{\down\,\,y} \ar[r]^{\widehat a} & A \ar[dl]^m \\
                                                     & X                                           &               }
\eeq
\end{enumerate} 
Thus, we have mappings $\,\ov X(x,y)\times y\ov m \to x\ov m$, and we now show how they can be enriched
as maps in $\S$.
(This enriching should not be confused with the $m_{a,b}:A(a,b)\to X(ma,mb)$ 
of Proposition~\ref{44a}, associated to any map in $\C$.)
Following and generalizing the results of Section~\ref{enrich}, we define $[xm]$ as the pullback
\eq   \label{60}
\xymatrix@R=4pc@C=4pc{
[xm] \ar[r]^{q_{x,m}}\ar[d]_{p_{x,m}}  & A \ar[d]^m \\
x\bs X  \ar[r]^{\up\,\,x}              & X          }
\eeq
and $(xm):=\comp[xm]$, with reflection map $c_{x,m}:[xm]\to (xm)$.
Then we also have the pullbacks:
\eq   \label{61}
\xymatrix@R=4pc@C=4pc{
(xm)\ar[d]\ar[r]^{i_{x,m}} & [xm] \ar[r]^{q_{x,m}}\ar[d]^{p_{x,m}}  & A \ar[d]^m \\
1         \ar[r]^{i_x}     & x\bs X  \ar[r]^{\up\,\,x}              & X          }
\qq\qq
\xymatrix@R=4pc@C=4pc{
(xm) \ar[d]\ar[r]   & A \ar[d]^m \\
1    \ar[r]^x       & X          }
\eeq
where we can assume that $c_{x,m}\circ i_{x,m}$ is the identity of the set $(xm)$, and
which show that the elements of $(xm)$ correspond bijectively to those of $x\ov m$.

Now we define (see~(\ref{44}))
\( \mu_{x,y,m}:X(x,y)\times (y m) \to (x m) \)
as $c_{x,m}\circ\mu'$:
\eq   \label{62}
\xymatrix@R=2pc@C=2pc{
& (x,y)\times (ym) \ar[ddl]_l\ar[ddr]^r\ar@{..>}[d]^{\mu'} &   \\
& [xm] \ar[dl]\ar[dr]\ar[d]^{c_{x,m}} &  \\
x\bs X \ar[dr] & (xm) & A \ar[dl]      \\
& X &                              }
\eeq
where $\mu'$ is universally induced by the maps $l$ and $r$, the product projections 
followed by $p_{x,y}\circ e_{x,y}$ and $q_{y,m}\circ i_{y,m}$, respectively.
To prove that this is the desired enrichment,
for a fixed point $a\in A$ over $y$ we consider the map $[x,a]:[x,y]\to[xm]$:
\eq   \label{63}
\xymatrix@R=2pc@C=2pc{
[x,y] \ar@{..>}[dr]^{[x,a]}\ar[rr]^{q_{x,y}}\ar[dddr]_{p_{x,y}} & & X/y \ar[dddr]\ar[dr]^{\widehat a} \\
& [xm] \ar[rr]^{q_{x,m}}\ar[dd]^{p_{x,m}} & & A \ar[dd]^m     \\ \\
& x\bs X  \ar[rr]    & & X          }
\eeq
and define $(x,a):=\comp[x,a]:(x,y)\to(xm)$, getting:
\eq   \label{64}
\xymatrix@R=4pc@C=4pc{
(x,y) \ar@/_1.5pc/[dd]\ar[r]^{i_{x,y}}\ar@{..>}[d]^{(x,a)} & [x,y]\ar[r]^{c_{x,y}}\ar[d]^{[x,a]}   & (x,y) \ar@{..>}[d]^{(x,a)} \\
(xm) \ar[r]^{i_{x,m}}\ar[d]                                & [xm]\ar[r]^{c_{x,m}}\ar[d]^{p_{x,m}}  & (xm)                \\
1     \ar[r]^{i_x}                                         & x\bs X                                              }
\qq\qq   
\xymatrix@R=4pc@C=4pc{
(x,y) \ar[r]^{i_{x,y}}\ar[d]^{(x,a)} & [x,y]\ar[r]^{q_{x,y}}\ar[d]^{[x,a]} &  X/y \ar[d]^{\widehat a}  \\
(xm)  \ar[r]^{i_{x,m}}               & [xm]\ar[r]^{q_{x,m}}                &  A                     }   
\eeq
The right hand diagram above shows that $(x,a)$ acts on an arrow $\alpha:1\to X(x,y)$ 
in the same way as $\widehat a$ acts on the corresponding left arrow $q_{x,y}\circ i_{x,y}\circ\alpha:1\to X/y$.
Furthermore we obtain that 
\( (x,a) =\comp[x,a] = c_{x,m}\circ[x,a]\circ i_{x,y} = c_{x,m}\circ[x,a]\circ e_{x,y} = \mu\circ(\id\times a)\circ\la\id,!\ra \) 
Then, $\mu(\alpha,a)=(x,a)\circ\alpha$ corresponds to $(\ov m\alpha)a$, proving the thesis.

If $n:B\to X$ is another discrete fibration over $X$ and $\xi:m\to n$ a map over $X$, we have a corresponding 
morphism of discrete fibrations $\ov\xi:\ov m\to\ov n$ in $\CatX$, or equivalently a natural transformation. 
On the other hand, $\xi$ also induces maps $[x,\xi]:[xm]\to[xn]$:
\eq   \label{65}
\xymatrix@R=2pc@C=2pc{
[xm] \ar@{..>}[dr]^{[x,\xi]}\ar[rr]^{q_{x,m}}\ar[dddr]_{p_{x,m}} & & A \ar[dddr]^m\ar[dr]^\xi \\
& [xn] \ar[rr]^{q_{x,n}}\ar[dd]^{p_{x,n}} & & B \ar[dd]^n     \\ \\
& x\bs X  \ar[rr]    & & X          }
\eeq
and one can check that the maps $(x,\xi):=\comp[x,\xi]:(xm)\to (xn)$ define a morphism of 
left $\S$-modules on $X$, which enriches the natural transformation $\ov\xi$.

\section{Balanced category theory}
\label{balcat}

Balanced category theory is category theory developed in a balanced factorization category $\C$, 
playing the role of $\Cat$ with the comprehensive factorization systems. 
We here present just a few aspects of it; others still need to be analyzed.

\subsection{The underlying functors}

Summarizing the results of Section~\ref{bal}, we have two ``underlying" functors: the ``left" one
\( \ov{(-)}:\C\to\Cat \) and the ``right" one \( \ov{(-)}':\C\to\Cat \),
where $\ov X(x,y)=\CX(X/x,X/y)$ and $\ov X'(x,y)=\CX(x\bs X,y\bs X)$.
They are isomorphic up to the duality functor \( (-)\op:\Cat\to\Cat \):
\eq    \label{66}
\xymatrix@R=1.5pc@C=6pc{
                                        &  \Cat\ar[dd]^{(-)\op}   \\
\C\ar[ur]^{\ov{(-)}}\ar[dr]_{\ov{(-)}'} &                         \\
                                        &  \Cat                    }   
\eeq
By propositions~\ref{18b} and~\ref{24}, both the underlying functors preserve the terminal object and sets; 
the left one preserves final and initial points, discrete fibrations and opfibrations,
slices and coslices and left and right adjunctible maps, while the right one reverses them.
This apparent asymmetry is only the effect of our naming of the arrows in $\E$, $\E'$, $\M$ and $\M'$.

The duality of the underlying functors (corollaries~\ref{48} and~\ref{53}) is a consequence of the fact 
that for any object $X$ of $\C$ we have an $\S$-enriched weak category  
(that is, the associativity of the composition $\mu$ may not hold; see Section~\ref{weak}), and for any map $f:X\to Y$
an $\S$-enriched functor $f_{x,y}:X(x,y)\to Y(fx,fy)$ such that the following diagrams in $\Set$ commute:
\eq   \label{67}
\xymatrix@R=4pc@C=4pc{
{\ov X}'(z,y)\times{\ov X}'(y,x)\ar[d]^{\ov\mu'} & \C(1,X(x,y)\times X(y,z))\ar[r]^\sim\ar[l]_\sim\ar[d]^{\C(1,\mu)} 
& \ov X(x,y)\times\ov X(y,z)\ar[d]^{\ov\mu} \\
{\ov X}'(z,x) & \C(1,X(x,z))\ar[r]^\sim\ar[l]_\sim & \ov X(x,z)                                  }
\eeq   
\eq   \label{68}
\xymatrix@R=4pc@C=4pc{
{\ov X}'(y,x)\ar[d]^{\ov f_{y,x}'} & \C(1,X(x,y))\ar[r]^\sim\ar[l]_\sim\ar[d]^{\C(1,f_{x,y})} & \ov X(x,y)\ar[d]^{\ov f_{x,y}} \\
{\ov Y}'(fy,fx)                    & \C(1,Y(fx,fy))\ar[r]^\sim\ar[l]_\sim                     & \ov Y(x,y)                      }
\eeq 
Then, denoting by $\SCat$ the category of $\S$-enriched weak categories 
and $\S$-enriched functors, diagram~(\ref{66}) can be enriched as
\eq    \label{69}
\xymatrix@R=1.5pc@C=6pc{
                                          &  \SCat\ar[dd]^{(-)\op}   \\
\C\ar[ur]^{\ov{\ov{(-)}}}\ar[dr]_{\ov{\ov{(-)}}'} &                  \\
                                          &  \SCat                    }   
\eeq
where $\ov{\ov X}(x,y) = \ov{\ov X}'(y,x) = X(x,y)$, 
while $\ov{\ov f}_{x,y} = \ov{\ov f}'_{y,x} = f_{x,y}:X(x,y)\to X(fx,fy)$.

Furthermore, for any df $m:A\to X$, the corresponding presheaf on $\ov X$ can also be enriched 
as a left $\S$-module on the weak $\S$-category $X$,
and the functor $\MX\to\Set^{\ov X\op}$  can be enriched 
to a functor $\MX\to\Xmod(\S)$ (and dually for $\MoX\to\Set^{\ov X}$).

\subsection{The tensor functor}

Given $p,q\in\CX$, their {\bf tensor product} 
is the components set $\ten_X(p,q):=\comp(p\times_X q)$ of their product over $X$.
We so get a functor $\ten_X:\CX\times\CX\to\S$.
In particular, $\ten(x\bs X,X/y) = X(x,y)$ and, if $m\in\MX$, $\ten(x\bs X,m) = (xm)$.
If $m\in\MX$, $n\in\MoX$ and $\C = \Cat$, $\ten(m,n)$ gives the classical tensor product (coend) $n\otimes m$ 
of the corresponding set functors $m: X\op\to\Set$ and $n:X\to\Set$.

\subsection{The inversion law}

If $q:Q\to X$ is any map, and $n:D\to X$ is a dof over $X$, pulling back 
the $\EM$-factorization $q =\,\down q\circ e_q$ along $n$ we get:
\eq   \label{72}
\xymatrix@R=4pc@C=4pc{
n\times q \ar[r]^{n\,\times e_q}\ar[d]\ar[dr]|>>>>>>>>>{n'} & n\,\times\!\down q \ar[dl] \ar[rd]  &          \\
D\ar[rd]^n                                                  &              Q \ar[d]^q\ar[r]^{e_q} & N(q)\ar[dl]^{\down\,\,q}  \\
                                                            &                                   X &            }
\eeq
where $\times$ denotes the product in $\CX$.
Since $n'$ is a dof, by the rsl 
\[ n\times e_q: n\times q\,\to n\,\times\!\down q \] 
is a final map over $X$; 
thus its domain and codomain have the same reflection in $\MX$:
\begin{prop}    \label{73}
If $n\in\MoX$ and $q\in\CX$,
\( \down_X(n\times q) \iso\, \down_X(n\,\times\!\down_X q) \); 
dually, if $m\in\MX$,
\( \up_X(m\times q) \iso\, \up_X(m\,\times\!\up_X q) \)
\end{prop}
\epf
Similarly, given two maps $p:P\to X$ and $q:Q\to X$, we have the following
symmetrical diagram:
\eq   \label{74}
\xymatrix@R=4pc@C=4pc{
                   & \up p\times q \ar[r]\ar[d]^e                          & Q \ar[d]^{e_q}          \\
p\,\times\!\!\down q \ar[r]^i\ar[d] & \up p\,\times\!\down q \ar[r]\ar[d]  & N(q) \ar[d]^{\down\,\,q} \\
P\ar[r]^{i_p}          & N(p) \ar[r]^{\up\,\,p}                            & X          }
\eeq
Since objects linked by initial or final maps have the same reflection in $\Mu = \S$,
we get the following ``inversion law":
\eq   \label{75}
\ten(\up p, q) \iso \ten(\up p,\down q) \iso \ten(p,\down q)
\eeq

\subsection{The reflection formula}

If $p$ is a point $x:1\to X$, diagram~(\ref{74}) becomes:
\eq   \label{76}
\xymatrix@R=4pc@C=4pc{
                   & x\bs q \ar[r]\ar[d]^e                          & Q \ar[d]^{e_q}          \\
(x\down q) \ar[r]^i\ar[d]           & [x\down q] \ar[r]\ar[d]  & N(q) \ar[d]^{\down\,\,q} \\
1\ar[r]^{i_x}      & x\bs X \ar[r]^{\up\,\,x}                              & X          }
\eeq
and the inversion law becomes the ``reflection formula": 
\eq     \label{77}
(x\down q) = \ten(x\bs X,\down q) \iso \comp(x\bs q) = \ten(x\bs X,q) \in\S
\eeq 
giving the enriched value of the reflection $\down q\in\MX$ at $x$.
If the functor 
\eq     \label{77a}
\MX\to\S^{{\ov X}\op} \qq m \mapsto \ten(-\bs X,m)  
\eeq
reflects isomorphisms, then the formula $\ten(-\bs X,q)$ determines $\down q$ up to isomorphisms.
In particular, in this case, $q$ is final iff the inverse image $x\bs q = q^*(x\bs X)$ of any slice of $X$
is connected: $\comp(x\bs q)=1,\, x\in X$. 

If $\C=\Cat$, we find again the ``classical" formula~(\ref{3a}): 
\[ \down q \iso \comp(-\bs q) \] 
which gives the presheaf corresponding to the discrete fibration $\down q$.

\subsection{Dense maps}
\label{dense}

In Section~\ref{emcat}, we have already discussed dense maps, from a ``left-sided" perspective.
But density is clearly a ``balanced" concept.
Indeed, the following is among the characterizations of dense functors in $\Cat$
given in~\cite{dense}: 
\begin{prop}    \label{70}
If $f:X\to Y$ is dense, then for any arrow $\alpha$ of $Y$ the object $\alpha/\!/f$,
obtained by pulling back the arrow interval $[\alpha]$ along $f$, is connected.
If the functor~{\rm(\ref{77a})} reflects isomorphisms, the reverse implication also holds. 
\end{prop}
\pf
The pullback square factors through the two pullback rectangles below:
\eq   \label{71}
\xymatrix@R=2pc@C=2pc{
\alpha/\!/f \ar[rr]^{e'}\ar[dd]\ar[dr] & & [\alpha] \ar[dd]\ar[dr]^{m'} \\
                     & f/y \ar[rr]^<<<<<<<<<<e\ar[dl] & &  Y/y \ar[dl]  \\
X \ar[rr]^f                   & & Y                 }
\eeq
If $f$ is dense, the map $e$ is final by definition; $m'$ is a discrete opfibration, because it is the composite
of the component inclusion $[\alpha]\to[x,y]$ and the pullback projection $q_{x,y}:[x,y]\to Y/y$,
which are both discrete opfibrations (see diagrams~(\ref{37}) and~(\ref{42})).
Then, by the rsl also $e'$ is final, and the connectedness of $[\alpha]$ implies that of $\alpha/\!/f$.
For the second part, note that $[\alpha]$, having an initial point, is a coslice of $Y/y$: 
\eq    \label{71a}
[\alpha] \iso \alpha\bs(Y/y) 
\eeq
Thus, $\alpha/\!/f$ is the pullback $\alpha\bs e'$, and the fact that it is connected for any $\alpha$
implies, by the hypothesis, the finality of $e'$. 
\epf

\subsection{Exponentials and complements}

Till now, we have not considered any of the exponentiability properties of $\Cat$.
In this subsection, we assume that for any $X\in\C$, the discrete fibrations and opfibrations over $X$ are
exponentiable in $\CX$.
\begin{prop}    \label{78}
If $m\in\MX$ and $n\in\MoX$ then $m^n\in\MX$ and, dually, $n^m\in\MoX$.
\end{prop}
\pf
If $n$ is a dof over $X$ and $q:Q\to X$ is any object in $\CX$, by Proposition~\ref{73}, 
\[ n\times\eta_q : n\times q\,\to n\,\times\!\down q \] 
is a final map over $X$ (recall that $e_q:q\to\,\down q$ is the reflecton map of $q\in\CX$ in $\MX$,
that is the unit $\eta_q$ of the adjunction $\down_X\adj i_X:\MX\inc\CX$). 
Then it is orthogonal to any df over $X$, in particular to $m$.
In turn, this implies that any unit $\eta_q:q\to\,\down q$ is orthogonal to $m^n$,
since in general if $L\adj R:\A\to\B$, an arrow $f:B\to B'$ in $\B$ is orthogonal to the object $RA$
iff $Lf:LB\to LB'$ is orthogonal to $A$ in $\A$:
\eq   \label{79}
\xymatrix@R=4pc@C=4pc{
B \ar[r]^h\ar[d]_f           & A \ar[d]  \\
B' \ar[r]\ar@{..>}[ur]^{u_*} & 1          }
\qq\qq
\xymatrix@R=4pc@C=4pc{
LB \ar[r]^{h^*}\ar[d]_{Lf} & RA \ar[d]  \\
LB' \ar[r]\ar@{..>}[ur]^u  & 1          }
\eeq
Indeed, in one direction the right hand diagram above gives a unique $u$, whose transpose $u_*$
makes the left hand one commute: $u_*\circ f = (u\circ Lf)_* = (h^*)_* = h$.
The converse implication is equally straightforward. 

Then, taking $q=m^n$ we see that $\eta_{m^n}$ is a reflection map with a retraction, and
so $m^n\in\MX$ (see~\cite{pis3} or also~\cite{bor}).
\epf
Any $m\in\MX$ has a ``complement" 
\eq    \label{79a}
\ten(m,-)\adj\neg m : \S \to \CX \qq S\mapsto (S\times X)^m
\eeq
where $S\times X$ is the constant bifibration over $X$ with value $S$.
By Proposition~\ref{78}, $\neg m$ takes values in $\MoX$, giving a broad generalization of 
the fact that the (classical) complement of a lower-set of a poset $X$ is an upper-set. 
(These items are trated at length in~\cite{pis} and~\cite{pis2}.
There, we conversely used the exponential law to prove the inversion law in $\Cat$; however, 
the present one seems by far the most natural path to follow, also in the original case $\C=\Cat$.)

\subsection{Codiscrete and grupoidal objects}

An object $X\in\C$ is {\bf codiscrete} if all its points are both initial and final
(that is, all the (co)slice projections are isomorphisms).
An object $X\in\C$ is {\bf grupoidal} if all its (co)slices are codiscrete.
(Note that codiscrete objects and sets are grupoidal.)
Since these concepts are autodual and the underlying functors preserve (or dualize) slices,
it is immediate to see that they preserve both codiscrete and grupoidal objects 
(which in $\Cat$ have the usual meaning).

\subsection{Further examples of balanced factorization categories}
\label{ex}

\begin{itemize}
\item
If $\C$ is an $\EM$-category such that $\E$ is $\M$-stable
(or, equivalently, satisfying the Frobenius law; see~\cite{clem}), 
then $\C$ is a bfc with $\EM = \EMo$.
\item
In particular, for any lex category $\C$ we have the ``discrete" bfc $\,\Cd$ and the 
``codiscrete" bfc $\,\Cc$ on $\C$. In $\Cd$, $\E=\E'={\rm iso}\C$ and $\M=\M'={\rm ar}\C$.
Conversely in $\Cc$, $\E=\E'={\rm ar}\C$ and $\M=\M'={\rm iso}\C$.
In $\Cd$ all object are discrete: $\S=\C$, while $1$ is the only connected object.
Conversely, in $\Cc$ all objects are connected, while $1$ is the only set.
So, for any $X\in\Cd$ ($X\in\Cc$) $\id:X\to X$ ($!:X\to 1$) is the reflection map of $X$ in sets
and, for any $x:1\to X$, $X/x = x\bs X = x$ ($X/x = x\bs X = \id_X$). 
Then, for any $X\in\Cd$ ($X\in\Cc$), $\ov X$ is the (co)discrete category on $\C(1,X)$.

Furthermore, if $X\in\Cd$ and $x,y:1\to X$, then $X[x,y]$ is pointless if $x\neq y$, while $X[x,x]=1$.
So the same is true for $X(x,y)$.
On the other hand, if $X\in\Cc$, then $X[x,y]=X$ and $X(x,y)=1$, for any $x,y:1\to X$.
Indeed, for any arrow $u_{x,y}:x\to y$, the arrow interval $[u_{x,y}]$ is simply $x,y:1\to X$.

\item
On the category $\Pos$ of posets, we can consider the bfc in which $m\in\M$ ($\M'$) 
iff it is, up to isomorphisms, a lower-set (upper-set) inclusion,
and $e:P\to X$ is in $\E$ ($\E'$) iff it is a cofinal (coinitial) mapping in the classical sense:
for any $x\in X$, there is $a\in P$ such that $x\leq ea$ ($ea\leq x$).

The category of internal sets is $\S = 2$, and
the component functor $\Pos\to 2$ reduces to the ``non-empty" predicate.
For any poset $X$ and any point $x\in X$, the slice $X/x$ ($x\bs X$) 
is the principal lower-set (upper-set) generated by $x$.
As in $\Cat$, $\ov X$ is isomorphic to $X$ itself.
This example also shows that the underlying functor may not preserve final maps:
a poset $X$ is ``internally" connected iff it is not empty, while $\ov X \iso X$ may well 
be not connected in $\Cat$.

Given a map $p:P\to X$, a colimiting cone $\lambda:p\to x$ simply indicates that $x$ is the sup 
of the set of points $pa ,\, a\in P$. The colimit is absolute iff such a sup is in fact a maximum.
A map $f:X\to Y$ is adequate iff any $y\in Y$ is the sup of the $fx$ which are less than or equal to $y$.

Furthermore, given $x,y:1\to X$, $[x,y]\subseteq X$ has its usual meaning:
it is the interval of the points $z\in X$ such that $x\leq z\leq y$.
So, $X(x,y)$ is the truth value of the predicate $x\leq y$, and the enriching of $\ov X$ in internal sets
gives the usual identification of posets with categories enriched over $2$.

By restricting to discrete posets, we get the bfc associated to the epi-mono factorization system on $\Set$.

\item
Let $\Gph$ be the category of reflexive graphs, with the factorization systems $\EM$ 
and $\EMo$ generated by $t:1\to 2$ and $s:1\to 2$ respectively, where $2$ is the arrow graph
(see~\cite{pis} and~\cite{pis2}).
In this case $\ov X$ is the free category on $X$, and the underlying functor is the
free category functor $\Gph\to\Cat$.

Furthermore, if $X\in\Gph$ and $x,y:1\to X$, then $X[x,y]$ is the graph which has as objects
the pairs of consecutive paths $\la \alpha,\beta\ra$, with $\alpha:x\to z$ and $\beta:z\to y$, while 
there is an arrow $\la \alpha,\beta\ra \to \la \alpha',\beta'\ra$ over $a:z\to z'$ iff $\alpha'=a\alpha$ 
and $\beta=\beta'a$.
In $\Gph$, $\S=\Set$ and the components of $X[x,y]$ are the paths $x\to y$.
Thus, the enriching $X(x,y)$ of $\ov X(x,y)$ is in fact an isomorphism.

\end{itemize}

\section{Conclusions}
\label{conc}

There are at least three well-estabilished abstractions (or generalizations) of category theory:
enriched categories, internal categories and 2-category theory.
Each of them is best suited to enlighten certain of its aspects and to capture new instances; 
for example, monads and adjunctions (via the triangular identities) are surely 2-categorical concepts,
while enriched categories subsume many important structures and support quantifications.

Here, we have based our abstraction on final and initial functors and discrete fibrations and opfibrations, 
whose decisive relations are encoded in the concept of balanced factorization category; 
in this context, natural transformations and/or exponentials are no more basic notions.
Rather, we have seen how the universal exactness properties of $\C$ and those depending on the 
two factorization systems, complement each other in an harmonious way (in particular, pullbacks and 
the reflection in discrete (op)fibrations) to give an effective and natural tool for proving categorical facts.

While not so rich of significant instances different from $\Cat$ as other theories, we hope to have shown that 
balanced category theory offers a good perspective on several basic categorical concepts and properties,
helping to distinguish the ``trivial" ones (that is, those which depend only on the bfc structure of $\Cat$)
from those regarding peculiar aspects of $\Cat$ (such as colimits, lextensivity, power objects and the arrow object).
The latter have been partly considered in~\cite{pis3} and deserve further study.

Summarizing, balanced category theory is
\begin{itemize}
\item
{\em simple}: 
a bfc is a lex category with two reciprocally stable factorization systems
generating the same (internal) sets;
\item
{\em expressive}:
many categorical concepts can be naturally defined in any bfc;
\item
{\em effective}:
simple universal properties guide and almost ``force" the proving of categorical properties
relative to these concepts;
\item
{\em symmetrical} (or ``balanced"):
the category $\Cat$ in itself does not allow to distinguish an object from its dual;
this is fully reflected in balanced category theory;
\item
{\em self-founded}:
it is largely enriched on its own internal sets, providing in a sense its own foundation
(see for example~\cite{law03} and~\cite{law66}).   
\end{itemize}

\begin{refs}

\bibitem[Ad\'amek et al., 2001]{dense} J. Ad\'amek, R. El Bashir, M. Sobral, J. Velebil (2001), On Functors which are Lax Epimorphisms, 
{\em Theory and Appl. Cat.} {\bf 8}, 509-521. 

\bibitem[Borceux, 1994]{bor} F. Borceux (1994), {\em Handbook of Categorical Algebra 1 (Basic Category Theory)}, 
Encyclopedia of Mathematics and its applications, vol. 50, Cambridge University Press.

\bibitem[Clementino et al., 1996]{clem} M.M. Clementino, E. Giuli, W. Tholen, (1996), Topology in a Category: Compactness, 
{\em Portugal. Math.} {\bf 53}(4), 397--433.

\bibitem[Lawvere, 1966]{law66} F.W. Lawvere (1966), {\em The Category of Categories as a Foundation for Mathematics},
Proceedings of the Conference on Categorical Algebra, La Jolla, 1965, Springer, New York, 1-20.

\bibitem[Lawvere, 1994]{law94} F.W. Lawvere (1994), {\em Unity and Identity of Opposites in Calculus and Physics},
Proceedings of the ECCT Tours Conference.

\bibitem[Lawvere, 2003]{law03} F.W. Lawvere (2003), Foundations and Applications: Axiomatization and Education,
{\em Bull. Symb. Logic} {\bf 9}(2), 213-224.

\bibitem[Pisani, 2007a]{pis} C. Pisani (2007a), Components, Complements and the Reflection Formula, 
{\em Theory and Appl. Cat.} {\bf 19}, 19-40. 

\bibitem[Pisani, 2007b]{pis2} C. Pisani (2007b), Components, Complements and Reflection Formulas, preprint, math.\-CT/0701457. 

\bibitem[Pisani, 2007c]{pis3} C. Pisani (2007c), Categories of Categories, preprint, math.CT/0709.0837

\bibitem[Street \& Walters, 1973]{stw73} R. Street and R.F.C. Walters (1973), The Comprehensive Factorization of a Functor,
{\em Bull. Amer. Math. Soc.} {\bf 79}(2), 936-941.

\end{refs}

\end{document}